\documentclass[11pt]{amsart}
\usepackage{amssymb,amsmath,booktabs}
\usepackage[margin=1in]{geometry}
\usepackage{tikz}
\definecolor{holdsK}{RGB}{46,184,92}
\definecolor{unknownK}{gray}{0.80}
\definecolor{failsK}{RGB}{205,60,55}
\usepackage[colorlinks=true,urlcolor=blue,linkcolor=black,citecolor=black]{hyperref}
\hypersetup{pdftitle={Smale's mean value constant (1981) equals 4},pdfauthor={Dave Goel}}

\newtheorem{theorem}{Theorem}[section]
\newtheorem{proposition}[theorem]{Proposition}
\newtheorem{lemma}[theorem]{Lemma}
\theoremstyle{remark}
\newtheorem{remark}[theorem]{Remark}

\newcommand{\D}{\mathbb{D}}
\newcommand{\T}{\mathbb{T}}
\newcommand{\C}{\mathbb{C}}
\newcommand{\HH}{\mathbb{H}}
\newcommand{\dist}{\operatorname{dist}}

\title{Smale's mean value constant (1981) equals 4}
\author{Dave Goel}
\gdef\fnone{\raggedright Expanded version, with further results~\cite{Gexp}: \url{https://hexagonmath.org/2611.00283}\spacefactor=3000 }
\gdef\fntwo{\raggedright Theorem~\ref{thm:main} is proved in Lean~4 from Lean's three standard axioms alone, and checked with \texttt{lake~comparator}: \url{https://github.com/david44hello/paper_smale-mean-value-constant-equals-4/tree/v3.10/lean}\spacefactor=3000 }

\begin{document}

\maketitle\footnotetext[1]{\fnone}\footnotetext[2]{\fntwo}

\begin{center}
\footnotesize First version: November 6, 2026. This version (3.10): November 25, 2026.
\end{center}

\begin{abstract}
For a complex polynomial $p$ of degree $d\ge 2$ with $p(0)=0$ and $p'(0)\ne 0$, let $S(p)$ be the minimum of $|p(c)/(c\,p'(0))|$ over the critical points $c$ of $p$, and let $K(d)$ be the supremum of $S(p)$ over all such $p$ of degree $d$. In 1981 Smale proved that $K(d)\le 4$ for every $d$, and asked whether $4$ can be replaced by $1$. We prove that Smale's mean value constant $\sup_dK(d)$ equals $4$ and that $K(d)\to4$ as $d\to\infty$, combining a reduction due to Jatar and Ng with an explicit construction. Each $K(d)$ is known to be below $4$, so the supremum is not attained.\footnotemark[1] These results are verified in Lean.\footnotemark[2]
\end{abstract}

\medskip
{\footnotesize\leftskip=3pc\rightskip=3pc\noindent\textbf{AI Disclosure.} This work was done overwhelmingly by AI models, steered by the author in a continuous critical dialogue, variously with Google's Gemini, Anthropic's Claude and OpenAI's ChatGPT. The author takes full responsibility for the paper.\par}

\section{Introduction}

Let $p$ be a complex polynomial of degree $d\ge 2$ with $p(0)=0$ and $p'(0)\neq 0$. A \emph{critical point} of $p$ is a zero of $p'$; since $p'(0)\ne0$, no critical point is $0$. Compare, at each critical point $c$, the slope $p(c)/c$ of the chord from $0$ to $c$ with the slope $p'(0)$ at $0$, and put
\[
S(p):=\min_{p'(c)=0}\left|\frac{p(c)}{c\,p'(0)}\right| .
\]
Any point where the derivative does not vanish can be moved to $0$ by a translation, and subtracting a constant then gives $p(0)=0$, so the normalization loses nothing. Smale's mean value theorem \cite{Smale1981} says that $S(p)\le 4$ for every such $p$. For $d\ge 2$ let
\[
K(d):=\sup\{S(p):\deg p=d\},\qquad K^*:=\sup_{d\ge 2}K(d).
\]
Thus $K(d)$ is the best constant in Smale's theorem for degree $d$, and $K^*$, \emph{Smale's mean value constant}, is the best constant for all degrees. Jatar and Ng \cite{JatarNg} call it $M$. The polynomial $z+z^d$ has $S=1-1/d$, so $K^*\ge 1$, and Smale asked whether $K^*=1$ \cite{Smale1981,Smale2000}. Smale listed this question again, as the Mean Value Problem, in the addenda to the second version of \emph{Mathematical problems for the next century}~\cite[p.~21]{Smale1998}. Moreover $K(d)$ does not decrease with $d$ \cite{Crane2006}. The idea is that adding $\varepsilon z^{d+1}$ to $p$, with $\varepsilon$ small, moves the old critical points and the quotients $|p(c)/(c\,p'(0))|$ there only slightly, and creates one new critical point, of modulus of order $1/\varepsilon$, where this quotient is large. So $\lim_d K(d)$ exists and equals $K^*$.

\begin{theorem}\label{thm:main}
$\displaystyle K^*=\lim_{d\to\infty}K(d)=4$.
\end{theorem}

So for every number $a<4$ there is a polynomial with $|p(c)/(c\,p'(0))|>a$ at \emph{every} critical point~$c$.

The proof is not effective: for a given $a<4$ it gives no bound on the least $d$ with $K(d)>a$. The supremum is not attained: $K(d)\le 4^{(d-2)/(d-1)}<4$ \cite{BeardonMindaNg}, and $K(d)<4-2.263/\sqrt d$ for $d\ge 8$ \cite{Crane2007}.

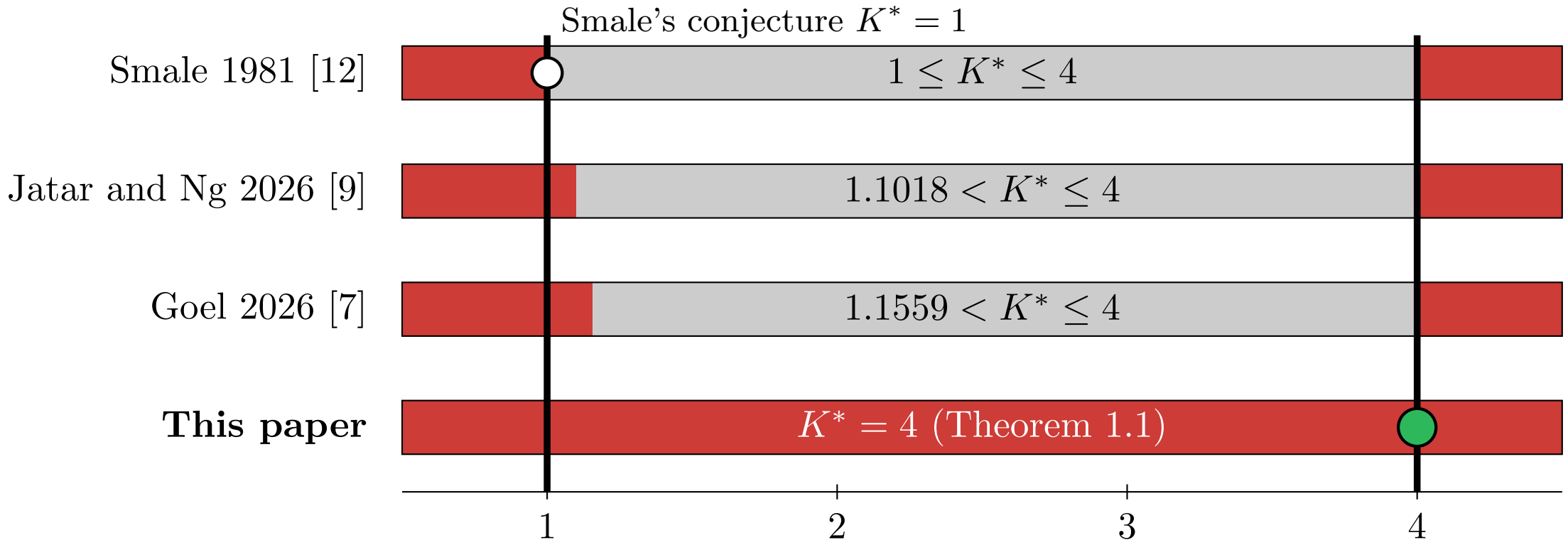
\begin{figure}[ht]
\centering
\begin{tikzpicture}[x=2.7cm,y=1cm,font=\small]
\fill[failsK] (0.5,3.3) rectangle (1,3.8);
\fill[unknownK] (1,3.3) rectangle (4,3.8);
\fill[failsK] (4,3.3) rectangle (4.5,3.8);
\draw (0.5,3.3) rectangle (4.5,3.8);
\node[left,inner sep=9pt] at (0.5,3.55) {Smale 1981 \cite{Smale1981}};
\node at (2.5,3.55) {$1\le K^*\le4$};
\node[above right, inner sep=2pt] at (1.02,3.8) {\footnotesize Smale's conjecture $K^*=1$};
\fill[failsK] (0.5,2.2) rectangle (1.1018586,2.7);
\fill[unknownK] (1.1018586,2.2) rectangle (4,2.7);
\fill[failsK] (4,2.2) rectangle (4.5,2.7);
\draw (0.5,2.2) rectangle (4.5,2.7);
\node[left,inner sep=9pt] at (0.5,2.45) {Jatar and Ng 2026 \cite{JatarNg}};
\node at (2.5,2.45) {$1.1018<K^*\le4$};
\fill[failsK] (0.5,1.1) rectangle (1.1559837,1.6);
\fill[unknownK] (1.1559837,1.1) rectangle (4,1.6);
\fill[failsK] (4,1.1) rectangle (4.5,1.6);
\draw (0.5,1.1) rectangle (4.5,1.6);
\node[left,inner sep=9pt] at (0.5,1.35) {Goel 2026 \cite{Goel1159}};
\node at (2.5,1.35) {$1.1559<K^*\le4$};
\fill[failsK] (0.5,0) rectangle (4.5,0.5);
\draw (0.5,0) rectangle (4.5,0.5);
\node[left,inner sep=9pt] at (0.5,0.25) {\textbf{This paper}};
\node[text=white] at (2.5,0.25) {\hypersetup{linkcolor=white}$K^*=4$ (Theorem~\ref{thm:main})};
\draw[line width=1.8pt] (1,-0.35) -- (1,3.9);
\draw[line width=1.8pt] (4,-0.35) -- (4,3.9);
\filldraw[fill=white, draw=black, thick] (1,3.55) circle (4pt);
\filldraw[fill=holdsK, draw=black, thick] (4,0.25) circle (5pt);
\draw (0.5,-0.35) -- (4.5,-0.35);
\foreach \k in {1,2,3,4} {\draw (\k,-0.28) -- (\k,-0.42) node[below] {$\k$};}
\end{tikzpicture}
\caption{What is known about $K^*$, to scale. In each row the red parts are ruled out, and the grey part is where $K^*$ may still lie; the black lines mark the interval $1\le K^*\le4$ of 1981. Smale proved $K^*\le4$ in 1981 and conjectured $K^*=1$ (white dot); that conjecture was disproved in 2026 \cite{LeanMVP, AdamczewskiColbrook}. Theorem~\ref{thm:main} shrinks the grey interval to the single point~$4$ (green dot).}
\label{fig:result}
\end{figure}

\subsection*{Notation}
$\D$ is the open unit disc and $\T$ the unit circle. A \emph{Jordan domain} is a bounded domain whose boundary is a Jordan curve, and it is \emph{analytic} if that curve is analytic. A function is \emph{univalent} if it is holomorphic and injective. Koebe's quarter theorem says that if $g$ is univalent on $\D$ and $g(0)=0$, then $g(\D)$ contains the disc of radius $|g'(0)|/4$ about $0$; the growth theorem gives the lower bound $|g(t)|\ge|g'(0)|\,|t|/(1+|t|)^2$ on $\D$. The \emph{Koebe function} $k(t)=t/(1+t)^2$ is univalent on $\D$ and maps it onto $\C\setminus[1/4,\infty)$. The boundary point of $k(\D)$ nearest $0$ is $1/4=\lim_{t\to1}k(t)$.

\subsection*{Earlier work}
The conjecture $K^*=1$ was disproved in 2026 \cite{LeanMVP, AdamczewskiColbrook}. Jatar and Ng then proved $K^*>1.1018$ \cite{JatarNg}, and the author proved $K^*>1.1559$ \cite{Goel1159} (Figure~\ref{fig:result}).

Our proof rests on a reduction due to Jatar and Ng \cite[Proposition~4.3, Theorem~4.4]{JatarNg} (Proposition~\ref{prop:bridge} below). The reduction replaces polynomials by functions with nonvanishing derivative near the closure of a domain; it goes back to the construction of \cite{LeanMVP, AdamczewskiColbrook} (see \cite[Proposition~6.4]{AdamczewskiColbrook}). Jatar and Ng show that $K^*$ equals the supremum $A$ of $\inf_{\partial\Omega}|f(z)/z|$ over such functions $f$ on analytic Jordan domains $\Omega\ni0$, normalized by $f(0)=0$ and $f'(0)=1$. In their notation our construction gives $A\ge4$, so $A=M=4$; this answers their question \cite[Section~9]{JatarNg} whether ``one can get a larger lower bound of $M$ through $A$''. Theorem~5.1 of \cite{JatarNg} shows that a different constant, the one in Smale's Problem~1A, equals~$4$; see Remark~\ref{rem:names}.

\subsection*{The idea of the proof}
By Proposition~\ref{prop:bridge}, it is enough to find, for each $a<4$, an analytic Jordan domain $\Omega\ni 0$ and a function $f$, holomorphic near $\overline\Omega$, with $f(0)=0$, $f'(0)=1$ and $f'\ne 0$ on $\overline\Omega$ such that $|f(z)/z|>a$ on $\partial\Omega$. For $a>1$ the function $f$ must have zeros other than $0$ on $\overline\Omega$, so it is not injective (Remark~\ref{rem:noninj}). We look for $f$ of the form $f=(W'(0)/F'(0))\,F\circ W^{-1}$, where $W$ is univalent on $\D$ with $W(0)=0$ and $\Omega=W(\rho\D)$ for some $\rho<1$ close to $1$. If $t_0$ is a point of the circle $\rho\T$ where $|F|$ is least, then the value of $|f(z)/z|$ at $z=W(t_0)$ is
\[
\frac{|F(t_0)|}{|F'(0)|}\cdot\frac{|W'(0)|}{|W(t_0)|}.
\]
The first factor is at most $\rho$, by Schwarz's lemma: the branch of $F^{-1}$ with value $0$ at $0$ extends to the disc $\{|v|<|F(t_0)|\}$ with values in $\rho\D$ (by analytic continuation of the inverse branch, using $F'\ne0$; see Remark~\ref{rem:why4}), so $|F(t_0)|\le\rho|F'(0)|$. The second is at most $(1+\rho)^2/\rho$, by Koebe's growth theorem. So this value is at most $(1+\rho)^2<4$. Our construction makes both bounds nearly sharp at the same point; this is possible because $F$ and $W$ can be chosen separately. The only link between them is a set of $n$ points of $\T$, the $n$-th roots $b_1,\dots,b_n$ of $-1$: near them $|F|$ dips, and there $|W|$ must be close to its smallest possible size.

\begin{figure}[ht]
\centering
\begin{tikzpicture}[scale=1]
\begin{scope}[xshift=-3.6cm]
\draw (0,0) circle (1.7);
\foreach \ang in {36,108,180,252,324} {
  \draw[line width=2.4pt] (\ang-5:1.7) arc (\ang-5:\ang+5:1.7);
}
\node at (1.72,1.25) [right] {\small $b_1$};
\node at (0,0.25) {$\D$};
\fill (0,0) circle (1.2pt) node[below] {$0$};
\node at (0,-2.25) {\small $|F|$ dips to about $1$ near the $b_j$};
\end{scope}
\draw[->] (-1.2,0) -- (-0.1,0) node[midway,above] {\small $\phi$};
\begin{scope}[xshift=2.2cm]
\fill[black!10] (0,0) circle (1.7);
\draw (0,0) circle (1.7);
\draw[line width=0.7pt] plot[smooth] coordinates {(1.698,0.057) (1.696,0.087) (1.691,0.134) (1.68,0.204) (1.652,0.31) (1.59,0.462) (1.458,0.663) (1.24,0.858)};
\draw[line width=0.7pt] plot[smooth] coordinates {(1.699,0.057) (1.697,0.087) (1.692,0.134) (1.681,0.205) (1.656,0.311) (1.598,0.466) (1.473,0.675) (1.229,0.912) (0.834,1.09) (0.359,1.096) (-0.044,0.927) (-0.121,0.87)};
\draw[line width=0.7pt] plot[smooth] coordinates {(1.699,0.057) (1.697,0.087) (1.693,0.134) (1.683,0.205) (1.659,0.312) (1.606,0.47) (1.489,0.687) (1.253,0.945) (0.854,1.163) (0.341,1.215) (-0.129,1.063) (-0.441,0.809) (-0.607,0.566) (-0.686,0.38) (-0.733,0.177) (-0.746,0.0)};
\draw[line width=0.7pt] plot[smooth] coordinates {(1.699,0.057) (1.697,0.087) (1.693,0.134) (1.684,0.205) (1.663,0.313) (1.614,0.473) (1.506,0.697) (1.282,0.977) (0.883,1.24) (0.332,1.349) (-0.212,1.226) (-0.599,0.958) (-0.814,0.681) (-0.918,0.46) (-0.981,0.216) (-0.998,-0.04) (-0.966,-0.294) (-0.889,-0.532) (-0.775,-0.744) (-0.617,-0.94) (-0.376,-1.137) (-0.325,-1.168)};
\draw[line width=0.7pt] plot[smooth] coordinates {(1.699,0.057) (1.697,0.087) (1.694,0.134) (1.686,0.205) (1.667,0.314) (1.623,0.476) (1.523,0.707) (1.313,1.007) (0.921,1.318) (0.338,1.498) (-0.291,1.423) (-0.775,1.15) (-1.057,0.833) (-1.198,0.567) (-1.285,0.268) (-1.308,-0.05) (-1.264,-0.364) (-1.159,-0.655) (-1.005,-0.907) (-0.798,-1.13) (-0.491,-1.338) (-0.091,-1.477) (0.356,-1.496) (0.777,-1.388) (1.112,-1.193) (1.343,-0.973)};
\fill (1.7,0) circle (1.4pt) node[right] {\small $1$};
\fill (0,0) circle (1.2pt) node[below] {$0$};
\node at (0,-2.25) {\small $\phi(\D)$, with its walls};
\end{scope}
\end{tikzpicture}
\caption{The construction for $n=5$, not to scale; the walls are drawn far apart. On the left is the unit disc; the thick arcs mark the neighbourhoods of the points $b_j$. On the right is $\phi(\D)$, the disc minus $5$ arcs of circles tangent to $\T$ at the point~$1$. The map $\phi$ sends points near each $b_j$ into a channel between neighbouring walls, or between the outermost wall and $\T$, close to the point $1$. The function $g$, close to the Koebe function, then sends the point $1$ to the boundary point nearest $0$.}
\label{fig:scheme}
\end{figure}
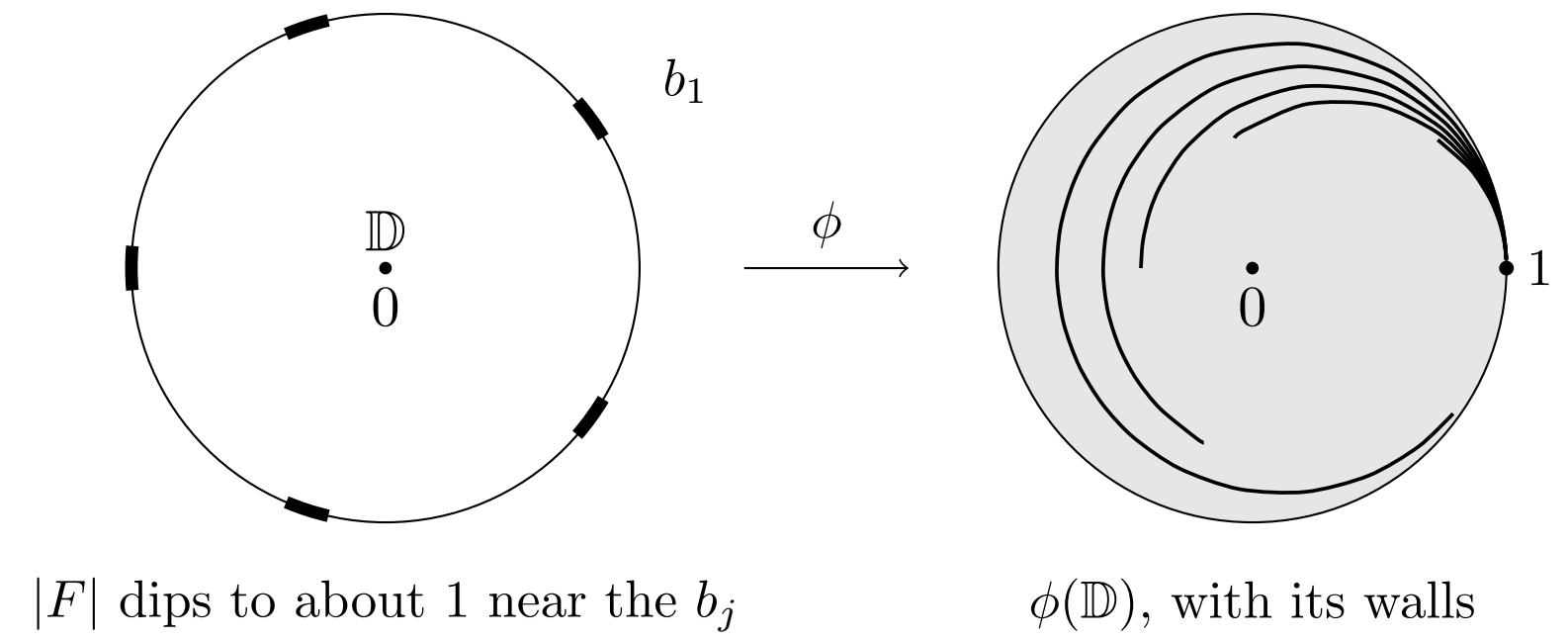

We take $W=g\circ\phi-g(\phi(0))$, where $g$ is a function close to the Koebe function and $\phi$ is an explicit univalent self-map of $\D$, with $\phi(0)$ close to $0$ and $\phi'(0)$ close to $1$, whose limit at \emph{each} of the points $b_j$ is the single point $1$. The image $\phi(\D)$ is the disc with $n$ thin walls removed, and the channels between the walls all lead to the point $1$. Since $g$ sends $1$ to the boundary point nearest $0$, the second factor is close to $4$ at each $b_j$. We construct $F$ with $F'(0)$ close to $1$ and with $|F|$ large near the whole circle, except near the points $b_j$, where $|F|$ drops to about $1$. Away from these points $|F|$ is so large that $|f(z)/z|$ is large. Our function $F$ is built from an entire function $E$ with $E'(\xi)=\exp(\xi^n)(1-\xi^n/c_{\sigma})$, where $c_{\sigma}$ is a large constant (Section~\ref{sec:F}). Figure~\ref{fig:scheme} shows the two pieces.

The main step is the following theorem, which we prove in Section~\ref{sec:proof}. For a function $g$ on $\D$ we write $g(1):=\lim_{t\to1,\,t\in\D}g(t)$, when this limit exists. We need the unrestricted limit, not only the radial or nontangential one, because $\phi(t)$ tends to $1$ tangentially to $\T$ as $t\to b_j$.

\begin{theorem}\label{thm:general}
Let $g$ be a bounded univalent function on $\D$ with $g(0)=0$ and $g'(0)=1$, and suppose that $g$ has a limit $g(1)$ at the point $1$. Then $K^*\ge 1/|g(1)|$.
\end{theorem}

The point $g(1)$ is not in $g(\D)$, since $g$ is open and injective, so it lies on $\partial g(\D)$; hence $|g(1)|\ge1/4$ by Koebe's theorem, and the lower bound $1/|g(1)|$ never exceeds $4$.

\section{From polynomials to pairs of functions on the disc}

The following proposition is due to Jatar and Ng \cite[Proposition~4.3]{JatarNg}. For completeness we include the proof.

\begin{proposition}\label{prop:bridge}
Let $\Omega$ be an analytic Jordan domain with $0\in\Omega$. Let $f$ be holomorphic on a neighbourhood of $\overline\Omega$, with $f(0)=0$, $f'(0)=1$ and $f'\neq 0$ on $\overline\Omega$. If $a>0$ and $|f(z)/z|>a$ on $\partial\Omega$, then $K^*>a$.
\end{proposition}

The proof uses the following form of Hilbert's lemniscate theorem \cite{Hilbert,NagyTotik}.

\begin{lemma}\label{lem:lemniscate}
In a closed topological annulus there is a nonconstant polynomial $q$ such that the whole set $\{|q|=1\}$ is an analytic Jordan curve in the interior of the annulus, separating its two boundary components. For such a $q$, $|q|<1$ inside this curve and $|q|>1$ outside it, and all zeros and all critical points of $q$ lie inside it.
\end{lemma}

\begin{proof}
We derive this single-curve form from the classical proof. Take an analytic Jordan curve in the interior of the annulus that separates its boundary components, and let $X$ be the closed region bounded by this curve. Let $G_X$ be the Green function of $\C\setminus X$ with pole at $\infty$; it has no critical point, since the complement of $X$ in the Riemann sphere is simply connected, so that $G_X$ is the logarithm of the modulus of a conformal map of $\C\setminus X$ onto $\{|w|>1\}$. Fix $\gamma>0$ so small that the set $\{0<G_X\le2\gamma\}$ lies in the interior of the annulus. Let $Q_j$ be a monic Fekete polynomial of degree $j$ of $X$, and put $q_j:=(e^{\gamma}\operatorname{cap}X)^{-j}Q_j$, where $\operatorname{cap}X$ is the capacity of $X$. Put $u_j:=j^{-1}\log|q_j|-(G_X-\gamma)=j^{-1}\log|Q_j|-\log\operatorname{cap}X-G_X$. The zeros of $Q_j$ (the Fekete points) lie on $\partial X$, so each $u_j$ is harmonic on $\C\setminus X$; it extends harmonically to $\infty$ with $u_j(\infty)=0$, since $Q_j$ is monic and $G_X(z)=\log|z|-\log\operatorname{cap}X+o(1)$. By the Bernstein--Walsh inequality, $u_j\le\varepsilon_j:=j^{-1}\log\|Q_j\|_X-\log\operatorname{cap}X$ on $\C\setminus X$, and $\varepsilon_j\to0$ \cite[Theorems~5.5.2, 5.5.4 and~5.5.7]{Ransford}. So $\varepsilon_j-u_j$ is a nonnegative harmonic function on $(\C\setminus X)\cup\{\infty\}$ with value $\varepsilon_j$ at $\infty$, and Harnack's inequality on the compact set $\{G_X\ge\gamma/4\}\cup\{\infty\}$ gives $0\le\varepsilon_j-u_j\le V_0\varepsilon_j$ there, with $V_0$ independent of $j$. Hence $u_j\to0$ uniformly on $\{G_X\ge\gamma/4\}$, including near $\infty$, and, by the interior estimates for harmonic functions, its first derivatives tend to $0$ uniformly on $\{\gamma/2\le G_X\le2\gamma\}$. Consider three regions, for large $j$.
\begin{itemize}
\item On the collar $\{\gamma/2\le G_X\le2\gamma\}$ the gradient of $G_X$ does not vanish. The gradient lines of $G_X$ make the collar a product of $\{G_X=\gamma\}$ with an interval. Along each of them, parametrized by arc length, $j^{-1}\log|q_j|$ increases (its derivative tends to $|\nabla G_X|>0$, uniformly on the compact collar) from about $-\gamma/2$ to about $\gamma$, so it vanishes exactly once. Hence $\{|q_j|=1\}$ meets the collar in a graph over $\{G_X=\gamma\}$ along the gradient lines, that is, in one Jordan curve close to $\{G_X=\gamma\}$; it is analytic, since $q_j\ne0$ and $\nabla\log|q_j|\ne0$ in the collar.
\item On $\{G_X\ge2\gamma\}$, up to $\infty$, we have $j^{-1}\log|q_j|\ge\gamma/2$, so $|q_j|>1$.
\item On $\{G_X\le\gamma/2\}$, which contains $X$, the Bernstein--Walsh inequality $|Q_j|\le\|Q_j\|_X\,e^{jG_X}$ \cite[Theorem~5.5.7(a)]{Ransford} and $\|Q_j\|_X^{1/j}\to\operatorname{cap}X$ \cite[Theorems~5.5.2 and~5.5.4]{Ransford} give $|q_j|<1$.
\end{itemize}
So $q:=q_j$ has the required property.

For such a $q$, $|q|<1$ inside this curve and $|q|>1$ outside it, by the maximum principle and connectedness. Since $|q|>1$ outside the curve, all $\deg q$ zeros of $q$ lie inside it; so $q$ maps the inside of the curve onto $\D$ as a proper map of degree $\deg q$, and by the Riemann--Hurwitz formula it has $\deg q-1$ critical points there, and these are all of them.
\end{proof}

\begin{proof}[Proof of Proposition~\ref{prop:bridge}]
\emph{Step 1: a polynomial close to $f$.} By continuity we can find $h>0$ and the closure $N$ of a Jordan domain, with $\overline\Omega\subset\operatorname{int}N$, such that $f$ is holomorphic near $N$, $f'\ne 0$ on $N$, and $|f(z)/z|>a+2h$ on $N\setminus\Omega$. Since $\C\setminus N$ is connected, Runge's theorem gives polynomials converging to $f$ uniformly near $N$, and by Cauchy's estimates their derivatives converge on $N$. Subtracting the constant term, we obtain a polynomial $p$ with $p(0)=0$, with $p'\neq0$ on $N$, with $|p(z)/z|>a+h$ on $N\setminus\Omega$, and with $|p'(0)|<(a+h)/(a+h/2)$.

\emph{Step 2: a lemniscate.} Apply Lemma~\ref{lem:lemniscate} to the closed annulus $N\setminus\Omega$. We obtain a polynomial $q$ whose lemniscate $\{|q|=1\}$ is a Jordan curve in the interior of this annulus. Then $\overline\Omega$ lies inside the curve, so $|q|<1$ on $\overline\Omega$. Moreover $zq'(z)\neq 0$ on the curve and outside it, because $0$ and all critical points of $q$ lie inside. Fix $\nu\in(0,1)$ so small that the compact set $A_\nu:=\{1-\nu\le|q|\le 1+\nu\}$ lies in $\operatorname{int}N\setminus\overline\Omega$, and so that $\{|q|<1-\nu\}$ contains $\overline\Omega$ and every critical point of $q$.

\emph{Step 3: the polynomials $H_m$.} Let $H_m:=p+zq^m$, a polynomial of degree $m\deg q+1$ once $m$ is large. Then $H_m(0)=0$, and
\[
H_m'=p'+q^{m-1}\,(q+mzq').
\]
We show that for all large $m$ every critical point of $H_m$ lies in $A_\nu$.

(i) On $\{|q|\le 1-\nu\}$, a compact subset of $N$, we have $m|q|^{m-1}\to 0$ uniformly, so $H_m'\to p'$ uniformly. Since $p'\neq 0$ on $N$, the polynomial $H_m$ has no critical point there.

(ii) On $\{|q|\ge 1+\nu\}$ the polynomial $zq'$ has no zero, and it has the same degree as $q$. Hence there are positive constants $V_1,V_2,V_3$ with $|zq'|\ge V_1(1+|z|)^{\deg q}$, with $|q|\ge V_2(1+|z|)^{\deg q}$ and with $|q/(zq')|\le V_3$ on this set. So $|q+mzq'|\ge (m/2)|zq'|$ for large $m$. Let $d_1:=\deg p'$, so that $|p'|\le V_4(1+|z|)^{d_1}$ on this set for some $V_4$, and let $m-1\ge d_1$. Since $|q|\ge1+\nu$ there, $|q|^{m-1}\ge(1+\nu)^{m-1-d_1}|q|^{d_1}$, and $|q|^{d_1}\ge V_2^{d_1}(1+|z|)^{d_1\deg q}\ge V_2^{d_1}(1+|z|)^{d_1}$. With $|q+mzq'|\ge(m/2)V_1(1+|z|)^{\deg q}\ge(m/2)V_1$ we get
\[
\bigl|q^{m-1}(q+mzq')\bigr|\ \ge\ (1+\nu)^{m-1-d_1}\,V_2^{d_1}(1+|z|)^{d_1}\cdot\frac{m}{2}\,V_1\ >\ V_4(1+|z|)^{d_1}\ \ge\ |p'|
\]
for large $m$, uniformly on this set. So there is no critical point there.

(iii) Therefore every critical point $c$ of $H_m$ lies in $A_\nu$. On $A_\nu$ we have $|zq'|\ge V_5>0$, $|q|\le 2$ and $|p'|\le V_6$, so $H_m'(c)=0$ gives
\[
|q(c)|^{m-1}=\frac{|p'(c)|}{|q(c)+mcq'(c)|}\le\frac{2V_6}{V_5 m}
\]
for large $m$, and hence $|q(c)|^m=O(1/m)$. Since $H_m(c)/c=p(c)/c+q(c)^m$ and $c\in N\setminus\Omega$, we obtain $|H_m(c)/c|\ge a+h-O(1/m)$, uniformly over the critical points. Finally $H_m'(0)=p'(0)+q(0)^m\to p'(0)$, because $|q(0)|<1$. Therefore
\[
\liminf_{m\to\infty}S(H_m)\ \ge\ \frac{a+h}{|p'(0)|}\ >\ a+\frac h2 ,
\]
and so $K^*>a$.
\end{proof}

We use Proposition~\ref{prop:bridge} in the following form.

\begin{lemma}\label{lem:pair}
Let $W$ be univalent on $\D$ with $W(0)=0$, and let $F$ be holomorphic on $\D$ with $F(0)=0$ and $F'\neq0$ on $\D$. If $0<\rho<1$, $s>0$, and
\[
|F(t)|\,|W'(0)|\ >\ s\,|F'(0)|\,|W(t)|\qquad\text{for all }t\text{ with }|t|=\rho,
\]
then $K^*>s$.
\end{lemma}

\begin{proof}
The set $\Omega:=W(\rho\D)$ is an analytic Jordan domain whose closure lies in $W(\D)$. The function $f:=(W'(0)/F'(0))\,F\circ W^{-1}$ is holomorphic on $W(\D)$, with $f(0)=0$, $f'(0)=1$ and $f'\neq0$. On $\partial\Omega$, with $z=W(t)$ and $|t|=\rho$, we have $|f(z)/z|=|W'(0)/F'(0)|\,|F(t)/W(t)|>s$. Proposition~\ref{prop:bridge} gives $K^*>s$.
\end{proof}

\section{The map $W$}

The map $\mu(t):=i(1+t)/(1-t)$ sends $\D$ conformally onto the upper half-plane $\HH$, with $\mu(0)=i$ and $\mu'(0)=2i$. Its inverse is $\mu^{-1}(w)=(w-i)/(w+i)$, so that $\mu^{-1}(w)-1=-2i/(w+i)$. Let $P\subset\T\setminus\{1\}$ be a finite nonempty set. For $b\in P$ the point $x_b:=\mu(b)$ is real. For $\lambda>0$ let
\[
\Gamma_\lambda(z):=z+\lambda\sum_{b\in P}\log(z-x_b)\quad(z\in\HH),\qquad \phi_\lambda:=\mu^{-1}\circ\Gamma_\lambda\circ\mu,
\]
with the principal branch of the logarithm.

\begin{lemma}\label{lem:comb}
\begin{enumerate}
\item[(a)] $\phi_\lambda$ is a univalent map of $\D$ into $\D$.
\item[(b)] For every $b\in P$, $\phi_\lambda(t)\to1$ as $t\to b$ in $\D$.
\item[(c)] As $\lambda\to0$, $\phi_\lambda(0)\to0$ and $\phi_\lambda'(0)\to1$.
\end{enumerate}
\end{lemma}

\begin{proof}
(a) We have $\operatorname{Im}\Gamma_\lambda(z)=\operatorname{Im}z+\lambda\sum_b\arg(z-x_b)>0$, so $\Gamma_\lambda$ maps $\HH$ into $\HH$. Further, $\Gamma_\lambda'=1+\lambda\sum_b1/(z-x_b)$ has negative imaginary part on $\HH$. Since $\HH$ is convex, the difference quotient of $\Gamma_\lambda$ between two points of $\HH$, being the mean of $\Gamma_\lambda'$ on the segment between them, has negative imaginary part, and so is not $0$. Hence $\Gamma_\lambda$ is injective.

(b) As $t\to b$, $z=\mu(t)\to x_b$ in $\HH$. Then $\operatorname{Re}\Gamma_\lambda(z)=\lambda\log|z-x_b|+(\text{a term bounded near }x_b)\to-\infty$, so $\Gamma_\lambda(z)\to\infty$, and $\phi_\lambda(t)-1=-2i/(\Gamma_\lambda(z)+i)\to0$. This holds however $t$ approaches $b$.

(c) As $\lambda\to0$, $\Gamma_\lambda(i)\to i$ and $\Gamma_\lambda'(i)\to1$. So $\phi_\lambda(0)=\mu^{-1}(\Gamma_\lambda(i))\to0$, and $\phi_\lambda'(0)=-4\Gamma_\lambda'(i)/(\Gamma_\lambda(i)+i)^2\to1$.
\end{proof}

Figure~\ref{fig:scheme} shows the shape of $\phi_\lambda(\D)$. On each interval of the real axis whose right end is a point $x_b$, including the unbounded interval $(-\infty,\min_bx_b)$, the imaginary part of $\Gamma_\lambda$ is constant and its real part runs from $-\infty$ up to a maximum and back. One checks that $\Gamma_\lambda(\HH)$ is $\HH$ minus horizontal rays going to the left, and $\phi_\lambda(\D)$ is $\D$ minus arcs of circles tangent to $\T$ at the point $1$: the walls of Figure~\ref{fig:scheme}. This description is not used in the proofs.

\begin{proposition}\label{prop:W}
Let $g$ be as in Theorem~\ref{thm:general}, and let $B:=\sup_\D|g|$. For every $\varepsilon\in(0,1)$ and every finite nonempty set $P\subset\T\setminus\{1\}$ we can find a univalent function $W$ on $\D$ with $W(0)=0$ and a number $\theta_0>0$ such that $|W'(0)|\ge1-\varepsilon$, $|W|\le2B$ on $\D$, and $|W(t)|<(1+\varepsilon)|g(1)|$ for every $t\in\D$ with $\dist(t,P)<\theta_0$.
\end{proposition}

Note that the bound $2B$ does not depend on $P$. This is what allows us, in Section~\ref{sec:proof}, to fix the size of $|F|$ before the number of points in $P$.

\begin{proof}
Let $W:=g\circ\phi_\lambda-g(\phi_\lambda(0))$. This is univalent on $\D$ by Lemma~\ref{lem:comb}(a), $W(0)=0$, and $|W|\le2B$. By Lemma~\ref{lem:comb}(c), as $\lambda\to0$ we have $g(\phi_\lambda(0))\to0$ and $W'(0)=g'(\phi_\lambda(0))\phi_\lambda'(0)\to1$. Fix $\lambda$ so small that $|W'(0)|\ge1-\varepsilon$ and $|g(\phi_\lambda(0))|<(\varepsilon/2)|g(1)|$. By Lemma~\ref{lem:comb}(b) and the limit of $g$ at $1$, for each $b\in P$ we have $|g(\phi_\lambda(t))-g(1)|<(\varepsilon/2)|g(1)|$ for all $t\in\D$ near $b$. Since $P$ is finite, this gives $\theta_0$.
\end{proof}

\section{The function $F$}\label{sec:F}

Throughout this section $n\ge8$ is an integer, $\sigma\ge100$, $c_{\sigma}:=2\sigma$, and $R\ge1$, with
\begin{equation}\label{eq:conditions}
200\,\sigma n\le e^{(\sigma-1)/2},\qquad R\le e^{(\sigma-1)/2}.
\end{equation}

\begin{proposition}\label{prop:F}
For all $\varepsilon,\theta>0$ we can find a function $F$, holomorphic on $\D$ with $F(0)=0$ and $F'\neq0$ on $\D$, and a number $r_0\in(0,1)$, such that $0<F'(0)\le(2\sigma)^{1/n}$ and, for every $t$ with $r_0<|t|<1$,
\[
|F(t)|\ge R,\qquad\text{or else both}\quad|F(t)|\ge\frac{n}{n+1}-\varepsilon\quad\text{and}\quad|t^n+1|<\theta .
\]
\end{proposition}

\begin{figure}[ht]
\centering
\begin{tikzpicture}[scale=1]
\begin{scope}[xshift=-3.3cm]
\fill[black!7] (0,0) circle (1.3);
\fill[white] (0,0) circle (0.95);
\draw (0,0) circle (1.3);
\draw[dashed] (0,0) circle (0.95);
\foreach \ang in {22.5,67.5,112.5,157.5,202.5,247.5,292.5,337.5} {
  \draw[line width=2.2pt] (\ang-7:1.13) arc (\ang-7:\ang+7:1.13);
}
\fill (0,0) circle (1.2pt) node[below] {\small $0$};
\node at (0,0.4) {\small $\D$};
\node at (0,-1.62) {\small $r_0<|t|<1$};
\end{scope}
\draw[->] (-1.75,0.15) -- (-0.35,0.15) node[midway,above] {\small $t\mapsto\eta(t^n)$};
\shade[left color=white, right color=black!12] (0,-1.8) rectangle (0.6,1.8);
\fill[black!12] (0.6,-1.8) rectangle (1.2,1.8);
\fill[black!12] (1.2,1.45) rectangle (6.4,1.8);
\fill[black!12] (1.2,-1.8) rectangle (6.4,-1.45);
\fill[black!35] (6.4,-1.8) rectangle (7,1.8);
\draw[dashed] (1.2,-1.45) rectangle (6.4,1.45);
\draw (0.5,1.8) -- (7,1.8) -- (7,-1.8) -- (0.5,-1.8);
\draw[dashed] (0,1.8) -- (0.5,1.8);
\draw[dashed] (0,-1.8) -- (0.5,-1.8);
\node at (3.8,0.35) {\small $\Pi_0$};
\fill (5.2,0) circle (1.2pt) node[below] {\small $0$};
\draw (1.2,-1.8) -- (1.2,-1.9) node[below] {\footnotesize $-C$};
\draw (6.4,-1.8) -- (6.4,-1.9) node[below left, inner sep=1pt] {\footnotesize $\sigma-1$};
\draw (7,-1.8) -- (7,-1.9) node[below right, inner sep=1pt] {\footnotesize $\sigma$};
\node[right, align=left] at (7.1,0) {\small $|\Phi|\ge R$\\[-1pt] \small (Step~3)};
\node[above, inner sep=2pt] at (3.8,1.8) {\small $|\Phi|\ge n/(n+1)$ (Step~4)};
\draw[thin] (0.8,-1.0) -- (0.35,-2.3);
\node[below left, inner sep=1pt] at (0.35,-2.3) {\small $|\Phi|\ge n/(n+1)-\varepsilon$ (Step~5)};
\end{tikzpicture}
\caption{The plan of Section~\ref{sec:F}, for $n=8$, not to scale. The map $t\mapsto\eta(t^n)$ sends $\D$ onto the half strip $\Pi$. Off the compact set $\Pi_0$, Steps~3--5 bound $|\Phi|$ from below (the band at the bottom edge is like the one at the top; in the far left, $\varepsilon$ comes from the choice of $C$ at the end of Step~6). For $r_0<|t|<1$ the point $\eta(t^n)$ lies outside $\Pi_0$, and it lies outside the dark band only when $t$ is close to an $n$-th root of $-1$ (thick arcs; Step~6).}
\label{fig:plan}
\end{figure}
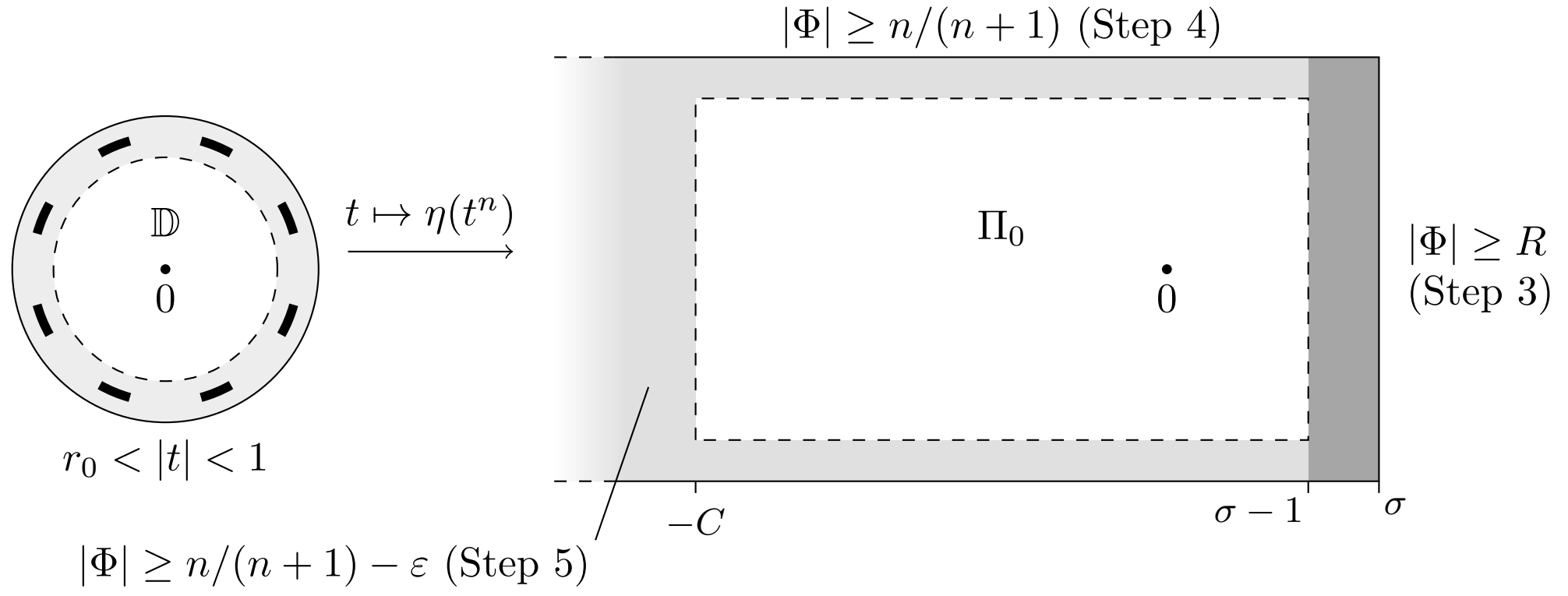

The plan of the proof is as follows (Figure~\ref{fig:plan}). We put $F:=E\circ\xi$, where $E$ is the entire function of Step~1 and $\xi(t)^n=\eta(t^n)$, with $\eta$ a conformal map of $\D$ onto the half strip $\Pi:=\{y:\operatorname{Re}y<\sigma,\ |\operatorname{Im}y|<C\}$, where $C$ is a large odd multiple of $\pi$, chosen last (Step~2; oddness is used in Step~4). Then $|F(t)|=|\Phi(y)|$ with $y=\eta(t^n)$, where $\Phi(y):=E(y^{1/n})$, and as $|t|\to1$ the point $y$ tends to the boundary of $\Pi$. In the band of width $1$ along the right edge $\operatorname{Re}y=\sigma$ we show $|\Phi|\ge R$ (Step~3). In the bands of width $\tfrac12$ along the two horizontal edges, and to the left of $\operatorname{Re}y=-C$, we show that $|\Phi|$ is at least about $n/(n+1)$ (Steps~4 and~5). The rest of $\Pi$ is a compact set, so for $|t|$ close to $1$ the point $y$ lies in one of these regions; it lies outside the right band only when $t^n$ is close to $-1$, that is, when $t$ is close to one of the $n$-th roots of $-1$ (Step~6). On $\Pi$ we have $\operatorname{Re}y<\sigma=c_{\sigma}/2$, so the factor $1-y/c_{\sigma}$ (Step~1) does not vanish there, and it is bounded below in terms of $|y|$ (Step~3): $|1-y/c_{\sigma}|\ge|y|/(6c_{\sigma})$. This lower bound is what lets Step~3 bound $|\Phi|$ from below along the whole right edge, at every height up to $C$, although $C$ is chosen last. The condition $n\ge8$ is used in Step~4.

\begin{proof}
\emph{Step 1: the function $E$.} Let
\[
E(\xi):=\int_0^\xi\exp(\zeta^n)\Bigl(1-\frac{\zeta^n}{c_{\sigma}}\Bigr)\,d\zeta .
\]
This is an entire function with $E'(\xi)=\exp(\xi^n)(1-\xi^n/c_{\sigma})$, and $E(\omega\xi)=\omega E(\xi)$ whenever $\omega^n=1$. Integrating by parts, with $E_0(\xi):=\int_0^\xi\exp(\zeta^n)\,d\zeta$ and $\kappa:=1/(c_{\sigma}n)$,
\begin{equation}\label{eq:closed}
E(\xi)=(1+\kappa)E_0(\xi)-\kappa\,\xi\exp(\xi^n).
\end{equation}
On $\C\setminus(-\infty,0]$, with principal branches and $u=\zeta^n$, put
\[
\Phi(y):=E(y^{1/n})=\int_0^yJ(u)\,du,\qquad J(u):=\frac1n\,e^u\chi(u),\qquad \chi(u):=u^{1/n-1}\Bigl(1-\frac{u}{c_{\sigma}}\Bigr).
\]
The function $\Phi$ itself jumps across the negative axis, but since $\Phi(\bar y)=\overline{\Phi(y)}$, its modulus $|\Phi|$ extends continuously to all of $\C$. Every $\xi$ with $\xi^n$ off the negative axis has the form $\omega\,(\xi^n)^{1/n}$ with $\omega^n=1$, so the rotation symmetry of $E$ and continuity give
\begin{equation}\label{eq:EPhi}
|E(\xi)|=|\Phi(\xi^n)|\qquad(\xi\in\C).
\end{equation}
So it suffices to bound $|\Phi|$ on the closed upper half-plane, with boundary values from above on the negative axis.

\emph{Step 2: the half strip and the definition of $F$.} Let $C\ge3\pi$ be an odd multiple of $\pi$, chosen last. Let $\beta:=\pi/(2C)$, $\alpha:=\sinh(\beta\sigma)$, and
\[
\Pi:=\{y:\operatorname{Re}y<\sigma,\ |\operatorname{Im}y|<C\}.
\]
The function $w(y):=\sinh(\beta(y-\sigma))$ maps $\Pi$ conformally onto the left half-plane: $y\mapsto\beta(y-\sigma)$ maps $\Pi$ onto the half strip $\{\operatorname{Re}z<0,\ |\operatorname{Im}z|<\pi/2\}$, which $\sinh$ maps conformally onto the left half-plane. So
\[
\Theta(y):=\frac{\alpha+w(y)}{\alpha-w(y)}
\]
maps $\Pi$ conformally onto $\D$, with $\Theta(0)=0$ since $w(0)=-\alpha$. Let $\eta:=\Theta^{-1}:\D\to\Pi$. Then $\eta(0)=0$, and from $w(\eta(t))=-\alpha(1-t)/(1+t)$ we get
\[
\eta'(0)=\frac{2\alpha}{\beta\cosh(\beta\sigma)}=\frac{2\tanh(\beta\sigma)}{\beta}\in(0,2\sigma].
\]
The quotient $\eta(t)/t$ has no zero on $\D$, so
\[
\xi(t):=t\left(\frac{\eta(t^n)}{t^n}\right)^{1/n},
\]
with the branch that is positive at $t=0$, is holomorphic on $\D$. It satisfies $\xi(t)^n=\eta(t^n)$ and $\xi(\omega t)=\omega\xi(t)$ for $\omega^n=1$. These two identities show that $\xi$ is injective: if $\xi(t_1)=\xi(t_2)$, then $\eta(t_1^n)=\eta(t_2^n)$, so $t_1^n=t_2^n$ and $t_2=\omega t_1$ with $\omega^n=1$; then $\xi(t_1)=\xi(t_2)=\omega\xi(t_1)$, so $\omega=1$ or $\xi(t_1)=0$, that is, $t_1=t_2=0$. Hence $\xi'\neq0$. Now let $F:=E\circ\xi$. Then $F(0)=0$, and
\[
F'(t)=\exp(\eta(t^n))\Bigl(1-\frac{\eta(t^n)}{c_{\sigma}}\Bigr)\xi'(t)\neq0,
\]
because $\operatorname{Re}\eta<\sigma<c_{\sigma}$. Also $F'(0)=\eta'(0)^{1/n}\le(2\sigma)^{1/n}$. By \eqref{eq:EPhi},
\begin{equation}\label{eq:Fboundary}
|F(t)|=|\Phi(\eta(t^n))|\qquad(t\in\D).
\end{equation}

\emph{Step 3: the right band, $\sigma-1\le\operatorname{Re}y<\sigma$: here $|\Phi(y)|\ge R$.} Let $y=x+i\tau$ with $\tau\ge0$ and $\sigma-1\le x<\sigma$, so that $|y|\ge x\ge99$. Since $\operatorname{Re}y\le c_{\sigma}/2$, we have $1+|y|/c_{\sigma}\le6|1-y/c_{\sigma}|$: if $|y|\le2c_{\sigma}$, then $|1-y/c_{\sigma}|\ge1-\operatorname{Re}y/c_{\sigma}\ge\tfrac12$ and $1+|y|/c_{\sigma}\le3$; if $|y|>2c_{\sigma}$, then $|1-y/c_{\sigma}|\ge|y|/c_{\sigma}-1\ge|y|/(2c_{\sigma})$ and $1+|y|/c_{\sigma}\le3|y|/(2c_{\sigma})$. Hence $|1-y/c_{\sigma}|\ge|y|/(6c_{\sigma})$. Let $T_0:=(1/n)e^x|\chi(y)|$, the modulus of the main term $(1/n)e^y\chi(y)$. Then
\[
T_0\ \ge\ \frac{e^x|y|^{1/n}}{6c_{\sigma}n}.
\]
To bound the rest of $\Phi(y)$, note that the function $\Phi(u)-(1/n)e^u\chi(u)$ has derivative $-(1/n)e^u\chi'(u)$, since $\Phi'=J=(1/n)e^u\chi$. Integrate this derivative along the horizontal segment from $y_0:=x/2+i\tau$ to $y$. For $u$ on this segment, $|y|/2\le|u|\le|y|$ and $|c_{\sigma}-u|\ge\sigma$. So $|\chi(u)|\le2(1+|y|/c_{\sigma})|y|^{1/n-1}\le12|\chi(y)|$ (the first step uses $|y|/2\le|u|\le|y|$, the second $1+|y|/c_{\sigma}\le6|1-y/c_{\sigma}|$), and $|\chi'(u)/\chi(u)|\le1/|u|+1/|c_{\sigma}-u|\le3/(\sigma-1)$. Therefore
\[
\Bigl|\Phi(y)-\frac1ne^y\chi(y)\Bigr|\ \le\ |\Phi(y_0)|+\frac1ne^{x/2}|\chi(y_0)|+\frac1n\int_{x/2}^{x}e^v\,\frac{36|\chi(y)|}{\sigma-1}\,dv .
\]
We bound the three terms on the right in turn.
\begin{itemize}
\item[(i)] Put $\xi_0:=y_0^{1/n}$, so that $\Phi(y_0)=E(\xi_0)$. On the segment from $0$ to $\xi_0$ we have $|\exp(\zeta^n)|\le e^{x/2}$, so $|E_0(\xi_0)|\le|\xi_0|e^{x/2}$; also $|\xi_0\exp(\xi_0^n)|=|\xi_0|e^{x/2}$. By \eqref{eq:closed}, $|y_0|\le|y|$ and the lower bound for $T_0$,
\[
|\Phi(y_0)|\le(1+2\kappa)|y_0|^{1/n}e^{x/2}\le2|y|^{1/n}e^{x/2}\le12c_{\sigma}n\,e^{-x/2}\,T_0\le\frac{24}{200}\,T_0 ,
\]
the last step because $12c_{\sigma}n\,e^{-x/2}=24\sigma n\,e^{-x/2}\le24\sigma n\,e^{-(\sigma-1)/2}\le24/200$ by \eqref{eq:conditions}.
\item[(ii)] Since $|\chi(y_0)|\le12|\chi(y)|$, the second term is at most $12e^{-x/2}T_0\le10^{-4}\,T_0$, as $x\ge99$.
\item[(iii)] Since $\int_{x/2}^xe^v\,dv\le e^x$, the third term is at most $36T_0/(\sigma-1)\le(36/99)\,T_0$.
\end{itemize}
Since $24/200+10^{-4}+36/99<1/2$, the total is at most $T_0/2$, and
\[
|\Phi(y)|\ \ge\ \frac{T_0}{2}\ \ge\ \frac{e^x}{24\sigma n}\ \ge\ \frac{e^{\sigma-1}}{24\sigma n}\ \ge\ \frac{200}{24}\,e^{(\sigma-1)/2}\ \ge\ R .
\]

\emph{Step 4: the band $C-\tfrac12\le\operatorname{Im}y<C$: here $|\Phi(y)|\ge n/(n+1)$, if $10c_{\sigma}\le C-\tfrac12$.} Fix $\tau$ in this range, and let $u=v+i\tau$ run along the horizontal line. Write
\[
J(u)=\frac{1}{nc_{\sigma}}\,e^v\,\bigl(-e^{i\tau}\bigr)\,u^{1/n}\Bigl(1-\frac{c_{\sigma}}{u}\Bigr).
\]
Since $C$ is an odd multiple of $\pi$, the factor $-e^{i\tau}$ has argument in $[-\tfrac12,0]$. The factor $e^{-i\pi/n}u^{1/n}$ has argument in $(-\pi/n,0)$. And $|c_{\sigma}/u|\le1/10$, so $1-c_{\sigma}/u$ has argument of modulus at most $0.11$. Since $\tfrac12+\pi/8+0.11<\pi/2$, with $n\ge8$ the total argument of $e^{-i\pi/n}J(u)$ has modulus below $\pi/2$, so $\operatorname{Re}(e^{-i\pi/n}J)>0$ along the whole line, and $\operatorname{Re}(e^{-i\pi/n}\Phi)$ increases along it as $v$ increases. As $v\to-\infty$, $\Phi(v+i\tau)-\Phi(v)\to0$, where $\Phi(v)$ is the boundary value from above. Indeed, on the vertical segment from $v$ to $v+i\tau$, of length less than $C$, we have $|J|\le(1/n)e^{v}(1+(|v|+C)/c_{\sigma})$ once $|v|\ge1$, because then $|u|^{1/n-1}\le1$ and $|1-u/c_{\sigma}|\le1+(|v|+C)/c_{\sigma}$. So the integral of $J$ over the segment is at most $C$ times this bound, which tends to $0$. On the negative axis $\Phi(-r^n)=e^{i\pi/n}\int_0^r\exp(-\zeta^n)(1+\zeta^n/c_{\sigma})\,d\zeta$. So the limit of $\Phi$ at the far left of the line is $e^{i\pi/n}\Lambda_n$, where
\[
\Lambda_n:=\int_0^\infty\exp(-\zeta^n)\Bigl(1+\frac{\zeta^n}{c_{\sigma}}\Bigr)\,d\zeta\ \ge\ \int_0^1(1-\zeta^n)\,d\zeta=\frac{n}{n+1}.
\]
Hence $|\Phi(y)|\ge\operatorname{Re}(e^{-i\pi/n}\Phi(y))\ge n/(n+1)$.

\emph{Step 5: the far left, $\operatorname{Re}y\le-C$, $0\le\operatorname{Im}y<C$: here $|\Phi(y)|\ge n/(n+1)-(C/n)e^{-C}(1+2C/c_{\sigma})$.} Let $y=x+i\tau$. On the negative axis, $|\Phi(x)|\ge\int_0^1(1-\zeta^n)\,d\zeta=n/(n+1)$, because $|x|^{1/n}\ge1$. On the vertical segment from $x$ to $y$, of length at most $C$, we have $1\le|u|\le|x|+C$, so $|J(u)|\le(1/n)e^{x}(1+(|x|+C)/c_{\sigma})\le(1/n)e^{-C}(1+2C/c_{\sigma})$, the last step because the middle expression decreases as $|x|$ grows. Integrating $J$ along this segment gives $|\Phi(y)-\Phi(x)|\le(C/n)e^{-C}(1+2C/c_{\sigma})$, and the bound follows.

\emph{Step 6: where $\eta$ goes, and the conclusion.} Let
\[
\Pi_0:=\{y:-C\le\operatorname{Re}y\le\sigma-1,\ |\operatorname{Im}y|\le C-\tfrac12\},
\]
a compact subset of $\Pi$, and let $r_0:=(\max_{\Pi_0}|\Theta|)^{1/n}<1$. Let $r_0<|t|<1$ and $y:=\eta(t^n)$. Then $|\Theta(y)|=|t|^n>r_0^n=\max_{\Pi_0}|\Theta|$, so $y\notin\Pi_0$. If $\operatorname{Re}y>\sigma-1$, Step 3 and \eqref{eq:Fboundary} give $|F(t)|\ge R$. Otherwise $\operatorname{Re}y<-C$ or $|\operatorname{Im}y|>C-\tfrac12$, and Steps 4 and 5 bound $|F(t)|$ from below. In this case $t^n$ is close to $-1$. Indeed $|\sinh(a_1+ia_2)|^2=\sinh^2a_1+\sin^2a_2$ shows that $|w(y)|\ge\sinh(\pi/2)$ if $\operatorname{Re}y<-C$, and $|w(y)|\ge\cos(\pi/(4C))\ge\tfrac12$ if $|\operatorname{Im}y|>C-\tfrac12$. Since $1+t^n=2\alpha/(\alpha-w)$ and $|\alpha-w|\ge|w|$ for $\operatorname{Re}w<0$, we get $|t^n+1|\le4\alpha$.

Given $\varepsilon,\theta>0$, choose $C$ so large that $10c_{\sigma}\le C-\tfrac12$, that $(C/n)e^{-C}(1+2C/c_{\sigma})\le\varepsilon$, and that $4\alpha=4\sinh(\pi\sigma/(2C))<\theta$. This proves the proposition.
\end{proof}

\section{Proof of the theorems}\label{sec:proof}

\begin{proof}[Proof of Theorem~\ref{thm:general}]
Let $B:=\sup_\D|g|$, so that $B\ge|g(1)|$, and let $R:=2B/|g(1)|\ge2$. Let $\varepsilon\in(0,1)$. For a large integer $\ell$ put $\sigma:=\ell$ and $n:=\ell^2$. Since $R$ does not depend on $\ell$, taking $\ell$ large we have $\sigma\ge100$, $n\ge8$, the conditions \eqref{eq:conditions}, the bound $(2\sigma)^{1/n}\le1+\varepsilon$ (by Bernoulli's inequality, $(1+\varepsilon)^{\ell^2}\ge1+\ell^2\varepsilon\ge2\ell$), and $n/(n+1)\ge1-\varepsilon/2$.

Let $P:=\{b_1,\dots,b_n\}$ be the set of $n$-th roots of $-1$; it does not contain $1$. Let $W$ and $\theta_0$ be as in Proposition~\ref{prop:W}. By compactness there is $\theta>0$ with $|t^n+1|\ge\theta$ for all $t\in\overline\D$ with $\dist(t,P)\ge\theta_0$. Let $F$ and $r_0$ be as in Proposition~\ref{prop:F}, for $\varepsilon/2$ and this $\theta$, and fix $\rho\in(r_0,1)$. Let $|t|=\rho$.
\begin{itemize}
\item If $|t^n+1|<\theta$, then $\dist(t,P)<\theta_0$, so $|W(t)|<(1+\varepsilon)|g(1)|$. Also $|F(t)|\ge1-\varepsilon$, because $R\ge2$ and $n/(n+1)-\varepsilon/2\ge1-\varepsilon$. So $|F(t)|/|W(t)|>(1-\varepsilon)/((1+\varepsilon)|g(1)|)$.
\item If $|t^n+1|\ge\theta$, then $|F(t)|\ge R$ by Proposition~\ref{prop:F}, since its other case needs $|t^n+1|<\theta$. Also $|W(t)|\le2B$, so $|F(t)|/|W(t)|\ge R/(2B)=1/|g(1)|>(1-\varepsilon)/((1+\varepsilon)|g(1)|)$.
\end{itemize}
So Lemma~\ref{lem:pair} applies with $s=(1-\varepsilon)|W'(0)|/((1+\varepsilon)|g(1)|F'(0))$. Using $|W'(0)|\ge1-\varepsilon$ and $F'(0)\le1+\varepsilon$,
\[
K^*\ >\ \Bigl(\frac{1-\varepsilon}{1+\varepsilon}\Bigr)^{2}\frac{1}{|g(1)|}.
\]
Letting $\varepsilon\to0$ gives the theorem.
\end{proof}

\begin{proof}[Proof of Theorem~\ref{thm:main}]
Let $k$ be the Koebe function, and for $0<\delta<1$ let $g(t):=k(\delta t)/\delta=t/(1+\delta t)^2$. Then $g$ is univalent on $\D$ and continuous on $\overline\D$, $|g|\le(1-\delta)^{-2}$, $g(0)=0$, $g'(0)=1$, and $g(1)=(1+\delta)^{-2}$. Theorem~\ref{thm:general} gives $K^*\ge(1+\delta)^2$, and letting $\delta\to1$ gives $K^*\ge4$. Smale's theorem gives $K^*\le4$. So $K^*=4$, and $\lim_{d\to\infty}K(d)=K^*$ because $K(d)$ does not decrease with $d$ \cite{Crane2006}.
\end{proof}

\section{Remarks}

\begin{remark}[Why the method stops at 4]\label{rem:why4}
Let $F,W,\rho,s$ be as in Lemma~\ref{lem:pair}, normalized so that $F'(0)=W'(0)=1$. Then $s<4$. Let $r$ be the largest radius for which the branch $\psi$ of $F^{-1}$ at $0$ extends to the disc $\{|v|<r\}$ with values in $\rho\D$. The branch $\psi$ cannot map $\{|v|<r\}$ into a disc $\{|t|\le\rho'\}$ with $\rho'<\rho$. Indeed, otherwise, since $F'\ne0$, the inverse function theorem and compactness give one radius $r_*>0$ such that for every $z_0$ with $|z_0|\le\rho'$, $F$ has an inverse branch with values in $\rho\D$ on the disc of radius $r_*$ about $F(z_0)$, taking $F(z_0)$ to $z_0$. Taking $z_0=\psi(v_1)$ with $r-r_*/2<|v_1|<r$, these branches continue $\psi$ to discs that cover $\{r\le|v|<r+r_*/2\}$; two of these discs that meet do so in a convex set that meets $\{|v|<r\}$ (it contains the midpoint of their centres), where both continuations equal $\psi$. So $\psi$ would extend to $\{|v|<r+r_*/2\}$ with values in $\rho\D$, against the choice of $r$. Hence there are points $v_k$ with $|\psi(v_k)|\to\rho$; then $|v_k|\to r$, since $\psi$ maps each disc $\{|v|\le r'\}$ with $r'<r$ onto a compact subset of $\rho\D$. Passing to a subsequence, $\psi(v_k)$ tends to a point $t^*$ with $|t^*|=\rho$, and then $|F(t^*)|=\lim|v_k|=r$. In particular $r\ge\min_{|t|=\rho}|F(t)|$. Schwarz's lemma, applied to $\psi$, gives $r\le\rho$, so $\min_{|t|=\rho}|F(t)|\le\rho$: this is the bound on the first factor used in the introduction. The function $W\circ\psi$ is univalent on this disc, with derivative $1$ at $0$, and $W(t^*)$ is a boundary point of its image, so $|W(t^*)|\ge r/4$ by Koebe's theorem. The hypothesis gives $r=|F(t^*)|>s|W(t^*)|\ge sr/4$. This is Smale's argument, carried over to the pair $(F,W)$. Jatar and Ng use the same two bounds in \cite[Theorem~7.1, Remark~7.3]{JatarNg}.
\end{remark}

\begin{remark}[Why $f$ cannot be injective]\label{rem:noninj}
If $f$ in Proposition~\ref{prop:bridge} has no zero other than $0$ on $\overline\Omega$, for instance if it is univalent, then $z/f$ is holomorphic on $\overline\Omega$ with value $1$ at $0$, and the maximum principle gives $\min_{\partial\Omega}|f(z)/z|\le1$. So $\min_{\partial\Omega}|f(z)/z|>1$ requires $f$ to have zeros other than $0$ on $\overline\Omega$. The function $F$ of Section~\ref{sec:F} has such zeros $t$; the points $\eta(t^n)$ lie inside the compact set $\Pi_0$, since Steps~3--5 bound $|\Phi|$ from below off $\Pi_0$. By Jensen's formula for $F(t)/t$, the mean of $\log|F|$ over $\rho\T$ equals $\log(\rho|F'(0)|)+\sum_k\log(\rho/|t_k|)$, summed over the zeros $t_k\ne0$ of $F$ in $\rho\D$. For $\rho$ close to $1$, $|F|$ is large on most of $\rho\T$ and bounded below on the rest (Steps~3--6), so this mean is large, while $\log(\rho|F'(0)|)$ is near $0$. The zeros supply the difference.
\end{remark}

\begin{remark}[The constant of Smale's Problem 1A]\label{rem:names}
The constant of \cite[Theorem~5.1]{JatarNg}, called $K$ there, is the one in Smale's Problem 1A \cite{Smale1981}, in the form of \cite[Section~5]{JatarNg}:
\[
K_{\mathrm{1A}}:=\sup_{j\ge2}\ \sup_p\ \min_{p'(c)=0}\Bigl|\frac{p^{(j)}(0)}{j!}\Bigr|^{1/(j-1)}|p(c)|\qquad(p(0)=0,\ p'(0)=1),
\]
in which $|p(c)|$ is not divided by $|c|$. The proof of \cite[Theorem~5.1]{JatarNg} applies a variant of the reduction of Proposition~\ref{prop:bridge} \cite[Proposition~5.2]{JatarNg} to the inverse of the Koebe function. That argument does not give $K^*$, because for that function $|f(z)/z|$ is small at boundary points far from $0$.
\end{remark}

\section{Discussion}\label{sec:discussion}

Before this paper it was known that $1.1559<K^*\le4$ \cite{Goel1159,Smale1981}. Theorem~\ref{thm:main} gives $K^*=4$. We end with three questions.

\emph{How fast does $K(d)$ tend to $4$?} From above, $K(d)\le4^{(d-2)/(d-1)}$ \cite{BeardonMindaNg} and $K(d)<4-2.263/\sqrt d$ for $d\ge8$ \cite{Crane2007}. From below, Theorem~\ref{thm:main} gives only the limit: our proof is not effective.

\emph{Which is the least degree $d$ with $K(d)>1-1/d$?} In its full strength, Smale's conjecture is $K(d)=1-1/d$; this holds for $d\le4$ \cite[Section~1]{Crane2007}, and $K(d)>1$ for every $d\ge29{,}984$ \cite[Section~1 and Theorem~A.1]{AdamczewskiColbrook}. So this degree lies between $5$ and $29{,}984$.

\emph{Blaschke products.} In the notation of Jatar and Ng, $M=B$ \cite[Theorem~1.2]{JatarNg}, where $B$ is the supremum of constants $B_d$, $d\ge1$, defined with Blaschke products of degree at most $d$; they prove $B_1<\sqrt{7/5}$ and $B_d<4-529(\log 2)^2/(576\,d^2)$ \cite[Theorem~1.3]{JatarNg}. By Theorem~\ref{thm:main}, $B=4$, so $B_1<B$. This answers two of the three questions that \cite[Section~9]{JatarNg} calls open: whether $B_1=B$, and whether $M=4$. The third, whether $B_d$ is strictly increasing in $d$, remains.

\end{document}